\documentclass[a4paper,12pt]{article}

\usepackage{amsthm}
\usepackage{amsmath,amssymb,latexsym,amsfonts,mathrsfs}
\usepackage[dvipdfmx]{graphicx}

\usepackage{color}

\usepackage{multicol}

\usepackage{fancyhdr}
\usepackage[top=30truemm, bottom=30truemm, left=25truemm, right=25truemm]{geometry}

\usepackage{ulem}

\newcommand{\SCR}[1]{{\mathscr #1}}

\newcommand{\CAL}[1]{{\cal #1}
}

\newcommand{\D}[1]{{\mathscr D}( #1 )}

\theoremstyle{definition}
\newtheorem{Thm}{{\bf Theorem}}[section]

\newtheorem{Lem}[Thm]{{\bf Lemma}}
\newtheorem{Prop}[Thm]{{\bf Proposition}}
\newtheorem{Cor}[Thm]{{\bf Corollary}}
\newtheorem{Ass}[Thm]{{\bf Assumption}}

\newtheorem{Rem}[Thm]{{\bf Remark}}

\newcounter{Exami}

\newcommand{\Proof}[2][Proof]{
\begin{proof}[{\bf #1}]
#2
\end{proof}
}

\begin{document}

%%%%%%%%%%%%%%%%%%%%%%%%%%%%%%%%%%%%%%%%%%%%%%%%%%%%%%%%%%%%%%%%%%%%%%%%%%%%%%%%%%%%%%%%%%%%%%%%%%%%%%%%%%%%%%%%%%%%%%%%%%%
%author
%%%%%%%%%%%%%%%%%%%%%%%%%%%%%%%%%%%%%%%%%%%%%%%%%%%%%%%%%%%%%%%%%%%%%%%%%%%%%%%%%%%%%%%%%%%%%%%%%%%%%%%%%%%%%%%%%%%%%%%%%%%
\begin{flushleft}
{\bf \Large On Schr\"{o}dinger equation with square and inverse-square potentials
 } \\ \vspace{0.3cm} 
 \end{flushleft}
\begin{flushleft}
{\large Atsuhide ISHIDA}\\
{
Katsushika Division, Institute of Arts and Sciences, Tokyo University of Science, 6-3-1 Niijuku, Katsushika-ku, Tokyo 125-8585, Japan\\ 
Alfr\'ed R\'enyi Institute of Mathematics, Re\'altanoda utca 13-15, 1053 Budapest, Hungary\\
Email: aishida@rs.tus.ac.jp
}
\end{flushleft}
\begin{flushleft}
{\large Masaki KAWAMOTO}\\
{Research Institute for Interdisciplinary Science, Okayama University, 3-1-1, Tsushimanaka, Kita-ku, Okayama City, Okayama, 700-8530, Japan }\\
Email: {kawamoto.masaki@okayama-u.ac.jp}
\end{flushleft}
\begin{abstract}
In this paper, we study the linear and nonlinear Schr\"odinger equations with a time-decaying harmonic oscillator and inverse-square potential. This model retains a form of scale invariance, and using this property, we demonstrate the asymptotic completeness of wave operators and Strichartz estimates for linear propagators.
\end{abstract}

\begin{flushleft}


\end{flushleft}
%%%%%%%%%%%%%%%%%%%%%%%%%%%%%%%%%%%%%%%%%%%%%%%%%%%%%%%%%%%%%%%%%%%%%%%%%%%%%%%%%%%%%%%%%%%%%%%%%%%%%%%%%%%%%%%%%%%%%%%%%%%
%Intro
%%%%%%%%%%%%%%%%%%%%%%%%%%%%%%%%%%%%%%%%%%%%%%%%%%%%%%%%%%%%%%%%%%%%%%%%%%%%%%%%%%%%%%%%%%%%%%%%%%%%%%%%%%%%%%%%%%%%%%%%%%%
\section{Introduction}
We consider the Hamiltonian
\begin{align} \label{1}
H(t) := \frac{p^2}2 + \frac{\sigma (t) x^2}{2} - \Lambda |x|^{-2} 
\end{align} 
acting on $L^2({\bf R}^n )$, where $n \in {\bf N}$, $n \geq 3$, $x \in {\bf R}^n $, $p = -i \nabla $, $\Lambda \leq  (n-2)^2/8 - \lambda _0 $, for some $\lambda _0 >0$, and $x^2 =|x|^2 $ and $p^2 = |p|^2 = -\Delta $. The coefficient $\sigma (t)$ of the harmonic potential satisfies the following conditions: 
\begin{Ass} \label{A1}
Assume $\sigma (t) \in L^{\infty} ({\bf R} _t)$, $ \sigma (t)  > 0$ for all $t \in {\bf R}$, $\sigma (t)$ is differentiable in $t$ and that there exists $C >0$ such that 
\begin{align} \label{6}
\sup_{t \geq s } \left| \frac{\sigma (t) - \sigma (s)}{\sigma (s) (t-s)} \right| \leq C.
\end{align}  
\end{Ass}
When we consider the scattering theory and the Strichartz estimates, we further assume the following: 
\begin{Ass} \label{A2} Let $\zeta (t)$ be one of the fundamental solutions of $y ''(t) + \sigma (t) y (t) = 0$. Then there exists $r_0 >0$, such that for all $|t| \geq r_0$, $\zeta (t) > 0$ and 
\begin{align}\label{4}
\pm  \int_{ \pm r_0}^{ \pm \infty}  \zeta (t) ^{-2} dt = \infty , 
\end{align} 
where all double signs correspond. 
\end{Ass}
A typical example satisfying Assumptions \ref{A1} and \ref{A2} is $\sigma (t) = \sigma _1 t^{-2}$ for $|t| \geq r_0$ and $\sigma (t) = \mbox{const}$ for $|t| < r_0$, where $\sigma _1 \leq 1/4$. Then, for $|t| \geq r_0$, $\zeta (t)$ can be written as 
\begin{align*}
\zeta (t) = 
c |t|^{ \lambda}, \quad c \neq 0 , \ \lambda = \frac{1 - \sqrt{1-4 \sigma _1}}2 .
\end{align*} 

\begin{Rem}
Condition on $\sigma (t)$ arises from the issue of existence of a unitary propagator, while condition \eqref{4} is necessary to address scattering issues.
\end{Rem} 

When considering the issue related to scattering, the first concern is the well-posedness of the linear Schr\"{o}dinger equations $i \partial _t u = H(t) u$, $u(s, x) = u_0(s)$, i.e., the unique existence of the propagator $U(t,s)$ for $H(t)$. When we replace the singular potential $-\Lambda |x|^{-2}$ with $V(t,x)$ which is in the Kato class, the well-posedness of solutions is known in the following sense based on the result by Yajima \cite{Ya}. Let 
\begin{align*}
D_{\pm 1} := \{ u \in \SCR{S}' ({\bf R}^n) \, | \,L^{\pm 1/2} u \in L^2({\bf R}^n) \}, \quad L = p^2 + x^2 +1,
\end{align*}
where $\SCR{S}' ({\bf R}^n)$ denotes the tempered distribution space. Then for any $u_0 (s) \in D_{+1}$, there uniquely exists a solution $u(t,x)$ of $i \partial _t u = H(t) u$, $u(s, x) = u_0(s)$ such that $u \in C( [s, \infty) \, ; \, D_{+1} ) \cap C^1([s, \infty) \, ; \, D_{-1})$. 
However, because $- \Lambda |x| ^{-2} $ is not included in the Kato class, Yajima's result cannot be applied directly. Therefore, we focus on the time-independence of  $- \Lambda |x| ^{-2} $. Indeed, the key estimate to show the existence of $U(t,s)$ in Yajima's work (Theorem 3.2 of \cite{Ya}, see, also Theorem 5.2 and Remark 5.3 in Kato \cite{Kato}) is denoted by
\begin{align} \label{10/25-1}
L^{-1/2} \left( \hat{H}(t) - \hat{H}(s)  \right) L^{-1/2} \in \SCR{B}(L^2({\bf R}^n)), 
\end{align} 
where the condition on the Kato class is required, and $\hat{H}(t)$ indicates a self-adjoint extension of $H(t)$ (see, \eqref{3} and \eqref{2} in Section 2). However, if the condition
\begin{align} \label{10/25-2} 
L^{-1/2} \hat{H}(t) L^{-1/2} = L^{-1/2} \left(  \overline{p^2/2 + \sigma (t)x^2/2} - \Lambda \overline{ |x|^{-2} }  \right) L^{-1/2}
\end{align}  
holds, it is evident that \eqref{10/25-1} is satisfied. From this, we have the unique existence of $U(t,s)$ for $H(t)$. To justify \eqref{10/25-2}, we need to exclude the critical coefficient for Hardy's inequality, namely the case where $\lambda_0 =0$.

As a result, we obtain the following theorem:
\begin{Thm}\label{T1}
Suppose Assumption \ref{A1}. Then there uniquely exists a family of unitary operators $U(t,s)$ such that
\begin{itemize}
\item $U(t, \tau ) U(\tau,s ) = U(t,s)$, $t,s,\tau \in {\bf R}$.
\item $(t,s) \mapsto U(t,s)$ is strongly continuous on $\SCR{B}(D_{+1})$.
\item For any $\varphi \in D_{+1}$, $(t,s) \mapsto U(t,s) \varphi \in D_{-1} $ is of a class of $C^1$ and satisfies $i \partial _t U(t,s) \varphi = H(t)U(t,s) \varphi $ and $i \partial _s U(t,s) \varphi = - U(t,s) H(s) \varphi $
\end{itemize} 
\end{Thm}

By virtue of the well-defined nature of the unitary propagator, it becomes feasible to define the operator
\begin{align*}
\CAL{W}^{}(t,s) :=  U(s,t) U_0(t,s),
\end{align*} 
wherein $U_0(t,s)$ denotes the propagator for $\hat{H}_0(t) = \overline{ p^2/2 + \sigma (t) x^2 /2 }$. In this paper, our primary focus lies in examining the existence of the wave operators 
\begin{align*}
\CAL{W}^{\pm} := \mathrm{s-} \lim_{t \to \pm \infty}\CAL{W}(t,s)
\end{align*}
and the unitariness on $L^2({\bf R}^n )$. Because energy is not conserved for time-dependent systems, the completeness of wave operators is discussed for some specific models. One such model is the time-periodic model, where the periodicity conserves Floquet energy (see, for example, Yajima \cite{Ya2}, Korotyaev \cite{Ko} among others). Another scenario involves strong external forces; for instance, Adachi-Ishida \cite{AI} demonstrates how external forces acting as acceleration forcibly remove all bound states. Without such advantageous properties, demonstrating completeness becomes a significantly challenging task. In such cases, discussions often focus solely on the characterization of the range of the wave operators (see, Yafaev \cite{Yaf}, Kitada-Yajima \cite{KY}, Kawamoto\cite{Ka} among others). 

As for our model, the decomposition formula by \cite{Ko} can reduce $H(t)$ to $ ( p^2 /2 - \Lambda  |x| ^{-2} )/\zeta (t) ^2 $. While $ ( p^2 /2 - \Lambda  |x| ^{-2} )/\zeta (t) ^2 $ represents time-dependent energy, the propagator  conserves energy with  $\CAL{U}(t,s) \varphi(H) = \varphi (H) \CAL{U}(t,s) $, $H=p^2 /2 - \Lambda  |x| ^{-2} $, where $\varphi \in C_0^{\infty} ((0, \infty))$ and  
\begin{align*}
\CAL{U}(t,s): = e^{\left( -i \int _s^t \zeta (\tau) ^{-2} d \tau \right) H}. 
\end{align*}  
This is the rationale behind demonstrating asymptotic completeness for our model, despite lacking periodicity and not being subjected to strong forces.

The first result in this paper is the following:
\begin{Thm} \label{T3}
Suppose Assumptions \ref{A1} and \ref{A2}. Then the wave operators $\CAL{W}^{\pm}$ exist and are complete, namely 
\begin{align*}
\mathrm{Ran} \left( \CAL{W}^{\pm} \right) = L^2({\bf R}^n).
\end{align*}
\end{Thm}

As a consequence of the reduction scheme, one observes that the Strichartz estimates for $\mathcal{U}(t,s)$ become significant by leveraging the results from the Strichartz estimates for $e^{-itH}$, as demonstrated in \cite{BPSZ} and \cite{BoMi} (also see \cite{Mi}). Let $(q, r)$ denote an admissible pair satisfying \begin{align} \label{8}
\frac{1}{q} + \frac{n}{2r} = \frac{n}4, \quad q \geq 2,
\end{align}
and let $L^q_IL^r $ be a Bochner-Lebesgue space, namely $F \in L^q_IL^r$ if and only if $F \in \SCR{S}'({\bf R}^{1+n})$ and $\left\|F \right\|_{L^q_IL^r} < \infty$, where for $q<\infty$
\begin{align*}
\| F \|_{L^q_IL^r} &:= \left( \int_{ I} \left\| F(t, \cdot) \right\|^{q}_{L^r({\bf R}^n)} dt  \right)^{1/q}
\end{align*}
and for $q=\infty$
\begin{align*} 
\| F \|_{L^{\infty}_I,r} &:=  \sup_{t \in  I}  \left\| F(t, \cdot) \right\|_{L^r({\bf R}^n)}.
\end{align*}

Applying the reduction scheme using $\mathcal{U}(t,s)$, we can demonstrate the following theorem:
\begin{Thm} \label{T2}
Suppose Assumption \ref{A1} and \ref{A2}. Let $I_{r_0}= [r_0, \infty)$, $s\in I_{r_0}$ and $(\tilde{q}, \tilde{r})$ be admissible pairs. Then for all $\phi \in L^2({\bf R}^n)$ and $F \in L^{\tilde{q}'}_{ I_{r_0} } L^{\tilde{r}'}$, there exists a constant $C>0$ that is independent of $t,s$ such that 
\begin{gather} \label{10/9-3}
\left\| U( t, s) \phi \right\|_{L^q_{ I_{ r_0} }L^r } \leq  C \| \phi  \|_{L^2({\bf R}^n)}, \\ 
\left\| 
\int_{r_0}^{  \infty} U(t,\tau) F(\tau) d\tau 
\right\|  _{L^q_{I_{r_0}}L^r}  \leq C \left\| F \right\|_{L^{\tilde{q}'}_{I_{r_0}} L^{ \tilde{r}'}}  \label{10/9-4} 
\end{gather} 
hold, where $\cdot '$ denotes the H\"{o}lder conjugate. The case where $I_{r_0}=(-\infty,-r_0]$ can be shown similarly. 
\end{Thm}

\begin{Rem}
In our case, we consider only the scenario where $|t| \geq r_0$. In this region, the useful factorization formula of $U(t,s)$ can be applied. However, if we consider the Strichartz estimate over the entire time domain $t \in (-\infty, \infty)$, it becomes challenging to prove Theorem \ref{T2} because the factorization formula fails for $t \in (-r_0, r_0)$.\end{Rem}
\begin{Rem}
If the generalized potential $V(t,x)$ satisfies that $V(x):=\zeta(t)^{2} V(t, \zeta _2(t) x)  $ is independent of the time, and $V(x)$ satisfies Assumption 2.1 of \cite{BoMi}. Then we may obtain the same theorems with replacing the potential $- \Lambda |x|^{-2}$ to  $V(t,x)$. The typical example of such $V(t,x)$ is $V(t,x) = - \Lambda |x|^{-2} + c |\zeta (t) |^{-1 } |x|^{- 1}  \chi (|x|/ \zeta (t) ) $ with $0 < c < \lambda_0 (2 + \sup|\chi '|)^{-1}/2  $ where $\chi \in C^1 ({\bf R})$ satisfies $|\chi(s)|\leq|s|^{-1}$ and $|\chi' (s)|\leq|s|^{-2}$ on $|s| \geq 1$, see Example 2.2 of \cite{BoMi}.
\end{Rem}

In \cite{KaYo} and \cite{Ka2}, the time-weighted Bochner-Lebesgue space was introduced, and Strichartz estimates were considered in this space. If $\Lambda =0$, the result of Fujiwara \cite{Fu} allows obtaining the finite time Strichartz estimates, and furthermore, combining this with Theorem \ref{T2} yields the following corollary, which is an extension of the results in \cite{KaYo} and \cite{Ka2}:

\begin{Cor}
Suppose $\sigma (t) \in L^{\infty} ({\bf R})$, Assumption \ref{A2} and $\Lambda =0$. Then the following Strichartz estimates hold: For $s\in {\bf R}$ and for any admissible pairs $(q,r)$ and $(\tilde{q}, \tilde{r} )$, 
\begin{gather*}
\left\| U(t,s) \phi \right\|_{L^q_{\bf R}L^r} \leq  C \| \phi  \|_{L^2({\bf R}^n)}, \\ 
\left\| 
\int_{ 0}^{ \pm t} U(t,\tau) F(\tau) d\tau
\right\|  _{L^q_{\bf R}L^ r} \leq C \left\| F \right\|_{L^{\tilde{q}'}_{\bf R} L^{ \tilde{r}'}} 
\end{gather*} 
hold. Here, in the case of $\Lambda=0$, the same estimates also hold for the case of $n =1$ or $n=2$, excluding condition $(q,r,n) = (2,\infty,2)$.
\end{Cor}

As an application of Theorem \ref{T2}, we consider the nonlinear equation 
\begin{align}\label{10/9-1} 
\begin{cases} 
i \partial _t u(t,x)&= H(t) u(t,x) + \lambda |u(t,x)|^{\theta} u(t,x) \\  
u(s,x) &= u_s
\end{cases}
\end{align} 
with $t \geq s \geq r_0$, $\lambda \in {\bf R} \backslash \{ 0\} $ and $1 < \theta < 4/n$. By using the Strichartz estimate and the approach outlined in Tsutsumi \cite{Tsu}, we can immediately demonstrate the following corollary:
\begin{Cor}
Suppose Assumptions \ref{A1} and \ref{A2}. For any $u_s \in L^2({\bf R}^n)$, an admissible pair $(q,r)$ and $T>r_0$, there uniquely exists a solution to \eqref{10/9-1} such that $u \in C([r_0, T] ;\, L^2({\bf R}^n)) \cap L^q_{\mathrm{loc}}([r_0,T] ;\, L^r({\bf R}^n))$.
\end{Cor}

In this paper, we mathematically formalize the scattering of particles under the influence of a time-decaying harmonic oscillator and an inverse square potential. Furthermore, we proved asymptotic completeness in such time-dependent quantum mechanical systems. These results would be important not only from a mathematical perspective but also from a physical standpoint. Moreover, in recent years, there has been rapid development in nonlinear problems under the influence of inverse-square potentials (see Zhang-Zheng \cite{ZZ} for example), and the contribution of the Strichartz estimate has been significant to these advancements. Moreover, in the case of $\Lambda =0$, nonlinear problems are also being investigated by \cite{Ka3} and \cite{KaMi}, with significant contributions from the development of the Strichartz estimate. Hence, this study suggests the feasibility of considering problems involving time-decaying harmonic potentials in addition to these developments, thus contributing to advancements in future research.

\section{Reduction of problems related to existence of wave operators}
Before considering the scattering issues concerning equation \eqref{1}, it is necessary to examine the existence of the unique propagator for $H(t)$, denoted as $U(t,s)$, which satisfies the following conditions:
\begin{align*}
& i \partial _t U(t,s) = H(t) U(t,s), \quad i \partial _s U(t,s) = - U(t,s) H(s), \quad U(s,s) = \mathrm{Id}_{L^2}, \\ 
& U(t,\tau)U(\tau, s) = U(t,s), \quad U(t,s) D_{+1} \subset D_{+1}.
 \end{align*}   
In this paper, we adopt the following simplified notations: \\ 
$\bullet$ \, $\| \cdot \| _{L^{r} ({\bf R}^n) }= \| \cdot  \|_{L^r} $, $1 \leq r \leq \infty$. \\ 
$\bullet$ \, For $u,v \in L^2({\bf R}^n)$, ${ \displaystyle ( u, v) = (u,v)_{L^2({\bf R}^n)} = \int_{{\bf R}^n }  u(x) \overline{v(x) } dx } $. 

\subsection{Existence of propagators}
We first demonstrate the existence of wave propagators. Given that $\sigma(t) >0 $ for any fixed $t$, we define 
\begin{align*}
Q_{H_0 (t)} (u,v) := ( H _0 (t) u,  v ), \quad u,v \in C_0^{\infty}({\bf R}^n  ) , \quad H_0(t) = \frac{p^2}2 + \frac{\sigma (t) x^2 }{2}.
\end{align*} 
By the Friedrichs extension (see, for example, Theorem X.23. in \cite{RS2}), the closure of $Q_{H_0 (t)}$ corresponds to the quadratic form of the self-adjoint operators $\hat{H}_0(t)$, each of which represents a unique positive extension of $H_0 (t)$.  Consequently we obtain a unique self-adjoint operator $\hat{H}_0(t)$ with form domain $ \D{Q_{\hat{H}_0(t) }} =D_{+1}$. Next we consider the quadratic form 
\begin{align} \label{10/9-2}
Q_{H(t)}(u,v) =  \left(H(t) u ,v \right), \quad u,v \in  C_0^{\infty}({\bf R}^n  ). 
\end{align} 
By the Hardy inequality, we have
\begin{align*}
Q_{H(t)} (u,u) \leq \frac{\sigma (t)}2 \| x u \|_{L^2}^2  + C \| \nabla  u \|_{L^2 }^2 ,
\end{align*} 
and 
\begin{align*}
Q_{H(t)} (u,u) &\geq \frac12 \| \nabla u \|_{L^2}^2 + \frac{\sigma (t)}2 \| xu \|_{L^2}^2 - \Lambda \| |x|^{-1} u  \|_{L ^2} ^2
\\ & \geq  \left( \frac12 - \frac{4 \Lambda}{(n-2)^2 }  \right) \left\| \nabla u \right\|_{L^2}^2 + \frac{\sigma (t)}2 \left\| x u \right\|_{L^2}^2  
\\ & \geq  \frac{ 4 \lambda_0}{(n-2)^2} \left\| \nabla u \right\|_{L^2}^2 + \frac{\sigma (t)}2 \left\| x u \right\|_{L^2}^2 . 
\end{align*} 
This implies that, for any fixed $t$, there exists $0< a_0 < A_0$ such that  
\begin{align*}
a_0 Q_{H_0(t)}(u,u) \leq Q_{H(t)}(u,u) \leq A_0 Q_{H_0(t)} (u,u) .
\end{align*} 
Because $C_0^{\infty}({\bf R}^n )  $ is dense in $D_{+1}  \subset L^2({\bf R}^n)$, the closure of $Q_{H (t)}$ corresponds to the quadratic form of the self-adjoint operator $\hat{H} (t)$, which represents a unique positive extension of $H (t)$.  Consequently, we obtain a unique self-adjoint operator $\hat{H}(t)$ with form domain 
\begin{align*}
\D{Q_{\hat{H}(t) }}&= \left\{ u \in L^2({\bf R}^n)\, ;\, \mbox{There exists}\  \{ u_n \}_{n\in{\bf N}} \subset C_0^{\infty} ({\bf R}^n)\right.\\
&\qquad\left. \mbox{such that}\ Q_{H(t)}(u_n -u_m, u_n -u_m) \to 0 \ \mbox{and}\ u_n \to u \ \mbox{as}\ n,m \to \infty \right\}\\
&= \left\{ u \in L^2({\bf R}^n)\, ;\, \mbox{There exists}\  \{ u_n \}_{n\in{\bf N}} \subset C_0^{\infty} ({\bf R}^n)\right.\\
&\qquad\left. \mbox{such that}\ Q_{H_0(t)}(u_n -u_m, u_n -u_m) \to 0 \ \mbox{and}\ u_n \to u \ \mbox{as}\ n,m \to \infty \right\} \\
&= \D{Q_{\hat{H}_0(t)} } = D_{+1}.
\end{align*} 
Therefore, we have the following lemma whenever $\lambda _0 >0$:
\begin{Lem}
We have $\D{Q_{\hat{H}(t)}} = \D{Q_{\hat{H}_0(t) }} = D_{+1}$.
\end{Lem}
We now demonstrate the unique existence of the propagator. To achieve this, we use Theorem 3.2 by Yajima \cite{Ya}. For the application of this theorem, we define the Hilbert spaces $\CAL{Y} = D_{+1}$, $\CAL{X} =D_{-1}$, $\CAL{Y}_t $ with inner product $Q_{\hat{H}(t) } (u,v)$, and $\CAL{X}_t$ with inner product $Q_{\hat{H}(t)^{-1} }   (u,v)$. Then we aim to obtain the followings to employ Theorem 3.2 by Yajima \cite{Ya}:
\begin{itemize}
\item $ \{ \hat{H}(t) , t \in I  \} $ be a family of closed operators in $\CAL{X}$ with dense domain $\CAL{Y} \subset \D{\hat{H}(t)}$.
\item $\hat{H}(t)$ is bounded operator from $\CAL{Y}$ to $\CAL{X}$ and norm continuous with respect to $I$.
\item Hilbert spaces $\CAL{X}_t$ and $\CAL{Y}_t$ satisfies 
\begin{align}\label{3}
 \| u \|_{\CAL{Y} _t}  / \| u \|_{\CAL{Y}_s}   \leq e^{c |t-s|}, \quad \| u \|_{\CAL{X} _t}  / \| u \|_{\CAL{X}_s}   \leq e^{c |t-s|}
 ,  \quad u \neq 0.
\end{align}   
\item $\hat{H}(t)$ is self-adjoint in $\CAL{X}_t$ and a part $\tilde{H}(t)$ of $\hat{H}(t)$ in $\CAL{Y}_t$ is also self-adjoint in $\CAL{Y}_t$.
\end{itemize}
Here we only show \eqref{3} since the norm continuity of $\hat{H}(t)$ follows from the discussion on the proof of the Lemma \ref{L-2/5-1} below, and other conditions are justified by the positiveness of $\hat{H}(t)$, as discussed in the last part of Section 4 and Section 5 of \cite{Ya}. 

\begin{Lem}
We have $\D{Q_{\hat{H}(t)^{-1}} } = \D{Q_{\hat{H}_0(t) ^{-1}  }} = D_{-1}$.
\end{Lem} 
\Proof{
For fixed $ t_ 0$ and $|t| \leq t_0$, let $\sigma ' = \min \{ 1/2,  \sigma (t_0) /2\} $. Then by the uncertainty property for the harmonic oscillator, we have, for any $u \in C_0^{\infty} ({\bf R}^n)$, 
\begin{align*}
 Q_{H_0(t) } (u,u) \geq  \sigma ' \| u \|_{L^2}^2 ,
\end{align*}
which implies $\D{Q_{\hat{H}_0(t)^{-1}}}  \supset L^2({\bf R}^n)$. Moreover by applying Lemma 3.1 of \cite{Ya}, we get $\D{Q_{\hat{H}_0(t)^{-1}}} = D_{-1}$. Finally, by 
\begin{align*}
A_0 Q_{H_0(t)} (u,u) \geq Q_{H(t)}(u,u) \geq {a_0} Q_{H_0(t) } (u,u) \geq a_0 \sigma ' \| u \|_{L^2}^2 ,
\end{align*} 
we also have $ \D{Q_{\hat{H}(t)^{-1}}}  = \D{Q_{\hat{H}_0(t)^{-1}}} = D_{-1}$. 
}

\begin{Lem}\label{L-2/5-1}
We have
\begin{align*}
 \| u \|_{\CAL{Y} _t}  / \| u \|_{\CAL{Y}_s}   \leq e^{c |t-s|}, \quad \| u \|_{\CAL{X} _t}  / \| u \|_{\CAL{X}_s}   \leq e^{c |t-s|}
 ,  \quad u \neq 0.
\end{align*} 
\end{Lem} 
\Proof{
Through a straightforward calculation and \eqref{10/9-2}, we find that for $u \in C_0^{\infty} ({\bf R}^n  )$, 
\begin{align*}
\left| \left(
Q_{H(t)} (u,u) - Q_{H(s)} (u,u) \right) \right|
& = \left| \left( (H(t) - H(s)) u,u \right) \right| \\
& \leq \frac12  \left| \frac{\sigma (t) - \sigma (s)}{\sigma (s) (t-s)} \right|  \left| \left( (t-s)\sigma (s)|x| ^2u, u \right) \right| 
\\ & \leq   \left| \frac{\sigma (t) - \sigma (s)}{\sigma (s) (t-s)} \right| |t-s| \left|   Q_{H(s)} (u,u)  \right|,
\end{align*} 
holds. This, together with Assumption \ref{A1}, yields 
\begin{align}\label{2}
\left| 
Q_{H(t)} (u,u) 
\right|  \leq \left( 1 + C|t-s| \right) |  Q_{H(s)} (u,u)  | \leq e^{ C|t-s| } | Q_{H(s)} (u,u) |.
\end{align}
Now we demonstrate the second inequality in \eqref{3}. Because the domains of $\D{ Q_{\hat{H}(t)} }  $ and $ \D{ Q_{\hat{H}(s) } } $ are equivalent to $ D_{+1}$, by using the formal relation 
\begin{align*}
H(t)^{-1} - H(s) ^{-1} = H(s)^{-1} ( H(s) - H(t) ) H(t)^{-1}   = \frac{( \sigma (s) - \sigma (t))}2 H(s)^{-1} |x|^2  H(t)^{-1},
\end{align*} 
we obtain
\begin{align*}
\left| Q_{H(t)^{-1} } (u,u) \right|  \leq \left( 1 + C|t-s|  \right) \left| Q_{H(s)^{-1} } (u,u) \right|  \leq e^{C|t-s|}   \left| Q_{H(s)^{-1} } (u,u) \right|.
\end{align*} }
Therefore, Yajima's Theorem 3.2 holds for this pair $H(t)$ and $H_0(t)$. This implies that Theorem \ref{T1} holds.

\subsection{Reduction of propagators}
To demonstrate the completeness of wave operators using direct estimation of $U(t,0)$ is quite difficult, as we lack an energy conservation law for $H(t)$ and uniform boundedness in $t$ for $\| \nabla U(t,0) \phi \|$, where $\phi \in \SCR{S} (\mathbf{R}^n)$. Therefore, we opt to use the non-zero condition on $\zeta(t)$ for $|t| \geq r_0$ and reduce a new Hamiltonian, which is easier to handle. We first define $\tilde{H}$ and the quadratic form $Q_{\tilde{H}}(u,v)$ of $\tilde{H}$ as 
\begin{align*}
\tilde{H} := \frac{p^2}2 - \Lambda |x|^{-2} , \quad Q_{\tilde{H}} (u,v) = (\tilde{H} u,v), \quad u,v \in C_0^{\infty} ({\bf R}^n). 
\end{align*} 
Using the same arguments in the case of $H(t)$, we can obtain $H$ as a positive extension of $\tilde{H}$ with form domain $\D{Q_H} = H^{1} :=  \left\{ u \in L^2({\bf R}^n) \, | \, \| u \|_{H^{1}} := \left\| (1 +p^2)^{1/2} u \right\|_{L^2}  < \infty  \right\}$, and we can also obtain the unitary propagator $e^{-itH}$. Here we define for $t \geq s > r_0$, (or $-r_0 > s \geq t$), 
\begin{align*}
\tilde{U}(t,s) := \SCR{J}(t) \CAL{U}(t,s) \SCR{J} ^{\ast}   ( s),  \quad \SCR{J}(t) := e^{i \frac{\zeta '(t)}{2 \zeta (t)}x^2} e^{-i \left(  \log \zeta (t) \right) A}  , \quad  \CAL{U}(t,s) = e^{-i \left(  \int_{s}^t \zeta (\tau)^{-2} d\tau \right) H },
\end{align*}
where $A = (x \cdot p + p \cdot x) /2$. We remark that 
\begin{align*}
\left( e^{-i \left(  \log \zeta (t) \right) A} \phi \right) (x) = \frac{1}{(\zeta (t))^{n/2} } \phi (x/ \zeta (t))
\end{align*} 
and that $\SCR{J}(t) C_0^{\infty} ({\bf R} ) \subset  C_0^{\infty} ({\bf R} )$, which leads us to find on $ C_0^{\infty} ({\bf R} )$,
\begin{align*}
\SCR{J}(t) H \SCR{J}^{\ast} (t)&= e^{i \frac{\zeta '(t)}{2 \zeta (t)} x^2} \zeta (t)^2 \left( \frac{p^2}{2 } - \frac{\Lambda}{|x|^2 }    \right)
 e^{-i \frac{\zeta '(t)}{2 \zeta (t)} x^2} 
 \\ &= \frac{\zeta (t)^2}{2  }\left( p- \frac{\zeta '(t)}{\zeta (t)} x \right)^2  - {\Lambda}{\zeta (t) ^2} |x|^{-2} .
\end{align*} 
This result together with 
\begin{align*}
i \partial _t \SCR{J}(t) = \left( -\frac{\zeta ''(t)\zeta (t) - (\zeta '(t))^2 }{2\zeta (t)^2 } x^2 + \frac{\zeta '(t)}{2\zeta (t)} \left( x \cdot \left( p-  \frac{\zeta '(t)}{\zeta (t)} x\right)+ \left( p-  \frac{\zeta '(t)}{\zeta (t)} x\right) \cdot x   \right) \right) \SCR{J}(t)
\end{align*}
and 
\begin{align*}
i \partial _t \CAL{U}(t,s)=  \frac{H}{\zeta (t) ^2} \CAL{U}(t,s),
\end{align*} 
we find that $i \partial_t \tilde{U}(t,s) = H(t) \tilde{U}(t,s)$ holds on $ C_0^{\infty} ({\bf R} )$.
Unique existence of the propagator leads us to the following proposition:
\begin{Prop}\label{P1}
Let $ |t| \geq |s| \geq r_0 $. Then the factorization of the propagator
\begin{align*}
U(t,s) = \SCR{J}(t) e^{-i \left(  \int_{s}^t \zeta ( \tau )^{-2} d \tau  \right) H } \SCR{J} ^{\ast}   ( s),  \quad \SCR{J}(t) := e^{i \frac{\zeta '(t)}{2 \zeta (t)}x^2} e^{-i \left(  \log \zeta (t) \right) A}  
\end{align*} 
holds. Moreover $U(t,s) D_{+1} \subset D_{+1}$ holds by Theorem \ref{T1}.
\end{Prop}

\section{Proofs of the main results}
As a consequence of Theorem \ref{T1}, we can define the pair of propagators $U(0,t) U_0(t,0)$ as well as $U_0(0,t) U(t,0)$. Then by using the same argument as in Section 2 the pair of the propagators can be rewritten as
\begin{align*}
U(0, t)^{\ast} U_0(t,0) &= U(0, r_0 )  \tilde{U}( r_0, t ) \tilde{U}_0(t, r_0) U_0(r_0, 0 ) 
\\ &= U(0,r_0) \SCR{J} ^{} (r_0) e^{ i \left(  \int_{r_0}^t \zeta (s)^{-2} ds \right) H  } e^{ \left( -i \int_{r_0}^t \zeta (s)^{-2} ds \right) H_0  } \SCR{J}^{\ast} (r_0) U_0(r_0,0) ,
\end{align*} 
where $H_0 = p^2/2$. Hence, asymptotic completeness follows by showing 
\begin{align*}
\mathrm{Ran} \left(   \underset{t \to \pm \infty}{\mathrm{s-}\lim \ } e^{ i \left(  \int_{r_0}^t \zeta (s)^{-2} ds \right) H  } e^{ \left( -i \int_{r_0}^t \zeta (s)^{-2} ds \right) H_0  }  \right) = L^2({\bf R}^n), 
\end{align*} 
which is equivalent to demonstrating the existence of 
\begin{align} \label{5}
 \underset{t \to \pm \infty}{\mathrm{s-}\lim \ } e^{ i \left(  \int_{r_0}^t \zeta (s)^{-2} ds \right) H_0  } e^{ \left( -i \int_{r_0}^t \zeta (s)^{-2} ds \right) H} .
\end{align} 
Here, using the fact that $\sigma (H) = \sigma_{\mathrm{ac}} (H)$, the intertwining property, and a density argument, it suffices to demonstrate the existence of 
\begin{align*}
  \underset{t \to \pm \infty}{\mathrm{s-}\lim \ } \tilde{\varphi}(H_0) e^{ i \left(  \int_{r_0}^t \zeta (s)^{-2} ds \right) H_0  } e^{ \left( -i \int_{r_0}^t \zeta (s)^{-2} ds \right) H} \varphi(H)
\end{align*} 
for any $\varphi, \tilde{\varphi} \in C_0^{\infty} ((0, \infty))$ with $\varphi = \varphi \tilde{\varphi}$. To help show the above, we introduce the following lemma: 
\begin{Lem} [Burq-Planchon-Stalker-Tahvildar-Zadeh \cite{BPSZ}] \label{L2/28}
For any $u \in L^2({\bf R}^n)$, 
\begin{align*}
\int_{- \infty} ^{\infty} \left\| 
|x|^{-1}  e^{ -i \tau H }  u
\right\| _{L^2} ^2 {d\tau} \leq C \| u \|_{L^2}^2
\end{align*} 
holds.
\end{Lem}

\subsection{Proof of Theorem \ref{T3}}
To prove the existence of \eqref{5}, it suffices to prove
\begin{align*}
\underset{t \to \pm \infty}{\mathrm{s-}\lim \ }e^{itH_0}e^{-itH}
\end{align*}
exist by virtue of the condition \eqref{4}. This is the direct consequence of Lemma \ref{L2/28} and the smooth perturbation theory by Kato (Theorem XIII.24 \cite{RS3}).

\subsection{Proof ot Theorem \ref{T2}} 
From Section 2, we find that $U(t,s)$ can be simply factorized as $Je^{-i \tau H} J^{-1}$. The propagator $e^{-i \tau H}$ has been widely investigated, particularly regarding Strichartz estimates and their nonlinear applications, which have been considered in several papers. The aim of this section is to address similar issues under the additional effect of a time-decaying harmonic potential.

We first establish the global-in-time Strichartz estimates for $U(t,s)$, where \eqref{10/9-3} can be demonstrated through similar calculations as those showing \eqref{10/9-4}. Hence, we focus solely on \eqref{10/9-4}. Since $-\Lambda |x|^{-2} $ is in $ L^{n/2, \infty} ({\bf R}^n _x)  $
, one can employ Theorem 2.11 of \cite{BoMi}, which enable us to handle inhomogeneous Strichartz estimates for any admissible pairs $(q,r)$ and $(\tilde{q}, \tilde{r})$ (without the restriction $q,\tilde{q}>2$). Hence for all $F \in C_0^{\infty} (\mathbf{R}^{1+n})$, we have the following 
\begin{align} \label{K2/5-1}
\left\| \int_0^{a} e^{-i(t-s) H} F(s) ds \right\|_{L^q([0, \infty) ; L^r)} \leq C \left\| F \right\|_{L^{\tilde{q}'}([0, \infty) ; L^{\tilde{r}'})}. 
\end{align}
Then it follows from Proposition \ref{P1} that 
\begin{align*}
\left\| 
\int_{r_0}^{t} U(t,s) F(s) ds
\right\| _{L^q([r_0, \infty);L^r)} ^q
&= \left\| 
\int_{r_0}^{t} \SCR{J}(t) e^{ -i H \int_s^t \zeta (\tau)^{-2} d \tau   }  \SCR{J}^{\ast}(s)F(s) ds
\right\|  _{L^q([r_0, \infty);L^r)} ^q
\\ & = 
\int_{r_0}^{\infty}   \left\| 
\int_{r_0}^{t} \SCR{J}(t) e^{ -i H \int_s^t \zeta (\tau)^{-2} d \tau   }  \SCR{J}^{\ast}(s)F(s) ds
\right\| _{r} ^q dt \\ 
&=
\int_{r_0}^{\infty}   \zeta (t) ^{-qn ( 1/2 - 1/r)}   \left\|
 \int_{r_0}^{t}  e^{ -i H \int_s^t \zeta (\tau)^{-2} d \tau   }  \SCR{J}^{\ast}(s)F(s) ds
\right\| _{r} ^q dt .
\end{align*}
Noting that $f : t \mapsto \int_{r_0}^t \zeta (\tau) ^{-2} d \tau $ is a bijection from $[r_0, \infty) $ to $[0, \infty)$, we have its inverse $f^{-1}$. Then change of variables $t=f^{-1}(a)$ and $s=f^{-1}(b)$ yields the following identity:
\begin{align*}
& \int_0^{\infty}  \zeta (t) ^{-qn ( 1/2 - 1/r)}   \left\|
 \int_{0}^{a}  e^{ -i (a-b) H    }  \SCR{J}^{\ast}(f^{-1} (b))F(f^{-1} (b) ) \zeta (s)^2 db
\right\| _{r} ^q \zeta (t)^2 da 
\\ &= 
   \left\|
 \int_{0}^{a}  e^{ -i (a-b) H    }  \SCR{J}^{\ast}(f^{-1} (b))F(f^{-1} (b) ) \zeta (f^{-1} (b))^2 db
\right\| _{L^q([0, \infty );L^r)} ^q, 
\end{align*} 
where we use $qn(1/2 - 1/r) =2$ from \eqref{8}. Using %Christ-Kiselev's lemma \cite{CK} and 
the Strichartz estimate \eqref{K2/5-1} due to \cite{BoMi}, we find the last term in the above equation is not greater than 
\begin{align}\label{2/28-1}
C  \left\|
 \SCR{J}^{\ast}(f^{-1} (b))F(f^{-1} (b) ) \zeta (f^{-1} (b))^2 
\right\| _{L^{\tilde{q}'}([0,\infty);L^{\tilde{r}'})} ^q.
\end{align} 
Using the change of variable $b=f(s)$, we find \eqref{2/28-1} is equivalent to 
\begin{align*}
C\left(  \int_{r_0}^{\infty} \zeta (s)^{2 \tilde{q}'-2}  \left\| \SCR{J}^{\ast} (s) F(s)  \right\|_{\tilde{r}'}^{\tilde{q}'}  ds \right)^{q/\tilde{q}'}  &= C\left(  \int_{r_0}^{\infty} \zeta (s)^{2 \tilde{q}'-2+ \tilde{q}'n ( 1/2 - 1/\tilde{r}') }  \left\|  F(s)  \right\|_{\tilde{r}'}^{\tilde{q}'}  ds \right)^{q/\tilde{q}'} 
\\ &=  C\left(  \int_{r_0}^{\infty} \left\|  F(s)  \right\|_{\tilde{r}'}^{\tilde{q}'}  ds \right)^{q/\tilde{q}'} ,
\end{align*} 
where we employ 
\begin{align*}
2 \tilde{q}'-2+ \tilde{q}'n \left( \frac12 - \frac1{\tilde{r}'} \right) = 0.
\end{align*} 
Therefore we finally conclude that 
\begin{align*}
\left\| 
\int_{r_0}^{\infty} U(t,s) F(s) ds
\right\| _{L^q([r_0, \infty);L^r)}  \leq C \left(  \int_{r_0}^{\infty}  \left\|  F(s)  \right\|_{\tilde{r}'}^{\tilde{q}'}  ds \right)^{1/\tilde{q}'} = C \left\| F \right\| _{L^{\tilde{q}'}([r_0, \infty);L^{\tilde{r}'})},
\end{align*} 
which is the desired result. 

\bigskip
\noindent
\textbf{Acknowledgments.}
This work was supported by JSPS KAKENHI Grant Numbers JP20K03625, JP20K14328, JP21K03279 and JP24K06796.

\end{document}